\documentclass[11pt]{article}
\usepackage{amsmath,amsthm, amssymb, amsfonts,latexsym,amscd}
\usepackage{graphicx}
\usepackage{tikz}
\usetikzlibrary{shadings,patterns}
\newtheorem{theorem}{Theorem}[section]
\newtheorem{corollary}[theorem]{Corollary}
\newtheorem{lemma}[theorem]{Lemma}
\newtheorem{proposition}[theorem]{Proposition}
\newtheorem{definition}[theorem]{Definition}

\begin{document}
	\title{\bf A  conjecture on  different  central parts of  binary trees}
	\author{Dinesh Pandey\footnote{Supported by UGC Fellowship scheme (Sr. No. 2061641145), Government of India} \and Kamal Lochan Patra}
	\date{}
	\maketitle
\begin{abstract}
	Let $\Omega_n$ be the family of binary trees on $n$ vertices obtained by identifying the root of an rgood binary tree with a vertex of maximum eccentricity of a binary caterpillar. In the paper titled ``On different middle parts of a tree (The electronic journal of combinatorics, 25 (2018), no. 3, paper 3.17, 32 pp)", Smith et al. conjectured that  among all binary trees on $n$ vertices the pairwise distance between any two of center, centroid and subtree core is maximized by some member of the family $\Omega_n$. We first obtain the rooted binary tree which minimizes the number of root containing subtrees and then  prove this conjecture. We also obtain the binary trees which maximize these distances.\\



\end{abstract}

\section{Introduction}
All the graphs in this paper are simple, connected and undirected. A tree $T$ is a connected, acyclic graph. The vertex and edge set of $T$ is denoted by $V(T)$ and $E(T)$ respectively. A subtree of $T$ is a connected subgraph of $T.$ For $u,v \in V(T)$, the distance  $d_T(u,v)$ or simply $d(u,v)$, is the number of edges in the path joining $u$ and $v$. The distance between two subsets $X$ and $Y$ of $V(T)$ is denoted by $d_T(X,Y)$ and defined as $d_T(X,Y)=\min\{d(x,y):x \in X, y\in Y\}$. The \textit{eccentricity} of a vertex $v \in V(T)$ is denoted by $e(v)$ and defined as $e(v)= \max\{d_T(u,v): u \in V(T)\}$. The radius $rad(T)$ of T is defined as $rad(T) = \min\{e(v) : v \in V\}$ and the diameter $diam(T)$ of $T$ is defined as $diam(T) = \max\{e(v) : v \in V\}.$ It is clear that $diam(T) = \max\{d_T (u,v) : u,v \in V\}.$ By degree of a vertex $v$ in $T$, we mean the number of edges incident with  $v$ and we denote it as $\deg(v)$. A vertex $v$  is referred as a \textit{pendant vertex} if $\deg(v)=1.$

By specifying a vertex $r\in V(T)$, we call  $T$, a rooted tree with root $r$.  A \textit{binary tree} is a  tree in which every non-pendant vertex has degree $3$.  A \textit{rooted binary tree} is a tree in which the root has degree two and any other vertex is either a pendant vertex or a vertex of degree $3$. Note that the number of vertices in a binary tree is always even and every binary tree on $n$ vertices has $\frac{n+2}{2}$ pendant vertices.  The number of vertices in a rooted binary tree is always odd. 
 
Let $T$ be a rooted binary tree with root $r$. For a vertex $v \in V(T)$, the \textit{height} $ht(v)$ of $v$ is  defined by  $ht(v)=d(v,r)$. Let $P(r,v)$ denotes the path joining $r$ and $v$. For $u,v \in V(T)$, $v$ is called a successor of $u$ if $P(r,u) \subset P(r,v)$. If $v$ is a successor of $u$ and $u,v$ are adjacent, we call $v,$ a child of $u$ and $u,$  the parent of $v$. The height of the tree is denoted by $ht(T)$ and is defined as $ht(T)=\max\{ht(v): v\in V(T)\}$.  We call a vertex $v \in V(T)$ is at level $l$ if its height is $l$. 

We call a rooted binary tree to be ordered, if for $l \geq 1$, the vertices at level $l$ are put in a linear order such that if $u$ and $v$ are vertices at level $l+1$ with different parents then the orders of $u$ and $v$ at level $l+1$ are same as the order of their parents at level $l$.

\begin{definition}\cite{Sw}
	A rooted binary tree is called an rgood binary tree if
	\begin{itemize}
		\item[(i)] The heights of any two of it's pendant vertices differ by at most 1 and
		\item[(ii)] The tree can be ordered such that the parents of the pendant vertices at the highest level make a final segment in the ordering of the vertices at next to highest level.
	\end{itemize}
\end{definition}

A single vertex rooted binary tree is also rgood. All rgood binary trees on $n$ vertices are isomorphic and we denote it by $T^n_{rg}$.  A \textit{caterpillar} is a tree which has a path such that every vertex not on the path is adjacent to some vertex on the path.  A \textit{binary caterpillar} is a caterpillar which is also a binary tree. Note that a binary caterpillar on $n$ vertices has diameter $\frac{n}{2}$. 

In the literature, different kind of central parts of trees  are studied for different perspectives. In this paper we are interested in the following central parts of binary trees: Center, Centroid and Subtree core.  We recall the definitions of these central parts which can be traced back to different papers.  The {\it center} of a tree $T$ is the set of vertices having minimum eccentricity. We denote center of a tree $T$ by $C(T)$. An element of $C(T)$ is referred as a central vertex. 

For $v\in V(T)$, a branch at $v$ is a maximal subtree of $T$ containing $v$ as a pendant vertex. The {\it weight} of $v$ is the maximal number of edges in any branch of $T$ at $v.$ We denote weight of a vertex $v \in V(T)$ by $Wt_T(v)$ or simply by $Wt(v)$.  We use $Wt_B(v)$ for weight of the branch $B$ at  the vertex $v$. A vertex of minimal weight is called a centroid vertex of $T$ and the set of all centroid vertices is called the {\it centroid} of $T.$ We denote the centroid of $T$ by $C_d(T)$. The following result is due to Jordan \cite{J}. 
\begin{proposition}\label{prop:cnt}(\cite{H},Theorem 4.2,Theorem 4.3)
\begin{enumerate}
\item	The center of a tree  consists of either a single vertex or two adjacent vertices.
\item The centroid of a tree consists of either a single vertex or two adjacent vertices.
\end{enumerate}
\end{proposition}

It is straight forward that $C(T)$ intersects with center of every longest path in $T$. If $|C_d(T)| = 2$ and $C_d(T) = \{u,v\}$, then $n$ must be even and $Wt(u) = Wt(v) = \frac{n}{2}.$ Also, among the branches at $u$ (respectively, at $v$), the branch containing $v$ (respectively, $u$) has the maximum number of edges. 

For $v \in V(T)$, $f_T(v)$ is the number of subtrees of $T$ containing $v$. The \textit{subtree core} of $T$ is the set of vertices $v$ for which $f_T(v)$ is maximum. We denote the subtree core of a tree $T$ by $S_c(T)$. The following result is due to Sz\'ekely and Wang.

\begin{proposition}\label{S_c-01}(\cite{Sw},Theorem 9.1)
The subtree core of a tree consists of either a single vertex or two adjacent vertices. 
\end{proposition}
The subtree core is the most recently defined central part of a tree. To prove Proposition \ref{S_c-01}, the authors used the fact that the function $f_T$ is strictly concave in the following sense.
\begin{lemma}\label{LS-01}
 If $u, v, w$ are three vertices of a tree $T$ with $\{u, v\}, \{v, w\} \in E(T ),$ then $2f_T(v) - f_T(u)- f_T(w)>0.$
\end{lemma}

The concept of central parts in trees were started by Jordan(\cite{J})  in 1869 with the definitions of center and centroid. Later many researchers contributed to it by giving definitons of median(\cite{Z}), telephone center(\cite{M}), distance center(\cite{Ka}), characteristic set(\cite{Rm}) and subtree core(\cite{Sw})(also see \cite{G}). But the idea of studying the pairwise distances between them is comparatively  new. In last $15$ years, the distance between these middle parts in various class of trees have been studied by many researchers (see \cite{Afjk,Dp,P,Sswy} ). Here we consider the class of binary trees on $n$ vertices and the pairwise distances between the central parts center, centroid and subtree core.

\subsection{Crg tree $T_{rg}^{n,l}$}\label{num}

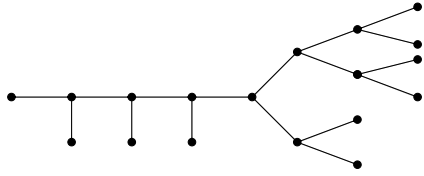
\begin{figure}[h]
\begin{center}
\begin{tikzpicture}
\filldraw (-.2,0)circle [radius=.5mm]--(0.6,0) circle [radius=.5mm]--(1.4,0)circle [radius=.5mm]--(2.2,0) circle [radius=.5mm]--(3,0) circle [radius=.5mm] ;
\filldraw (0.6,0)--(.6,-.6) circle [radius=.5mm];
\filldraw (1.4,0)--(1.4,-.6)  circle [radius=.5mm];
\filldraw (2.2,0)--(2.2,-.6)  circle [radius=.5mm];
\filldraw (3,0)--(3.6,.6)circle [radius=.5mm]--(4.4,.9)circle [radius=.5mm]--(5.2,1.2)circle [radius=.5mm];
\filldraw (4.4,.9)--(5.2,.7) circle [radius= .5mm];
\filldraw (3.6,.6)--(4.4,.3)circle [radius=.5mm]--(5.2,0)circle [radius=.5mm] ;
\filldraw (4.4,.3)circle [radius=.5mm]--(5.2,.5)circle [radius=.5mm] ;
\filldraw (3,0)--(3.6,-.6)circle [radius=.5mm]--(4.4,-.3)circle [radius=.5mm];
\filldraw (3.6,-.6)--(4.4,-.9)circle [radius=.5mm];
\end{tikzpicture}
\caption{The crg tree  $T_{rg}^{18,11}$} \label{fgr-00}
\end{center}
\end{figure}

Let $n\geq 4$ be a positve even integer and let $l\geq 3$ be a positve odd integer such that $l<n.$ Let $T_{rg}^{n,l}$ denote the tree on $n$ vertices which is obtained by identifying the root of $T_{rg}^l$ with a vertex of maximum eccentricity of a binary caterpillar tree on $n-l+1$ vertices (see figure \ref{fgr-00}). Such a  tree $T_{rg}^{n,l}$ is called a \textit{crg tree}. The binary caterpillar is the crg tree $T_{rg}^{n,3}$. 

We lable the vertices of a longest path of the caterpillar part of $T_{rg}^{n,l}$ by $1,2,\ldots, \frac{n-l+3}{2}=v$, where $v$ is the root of the rgood part of it. We denote by $\Omega_n$ the class of all  crg trees on $n$ vertices. Any binary tree on  $n\leq 8$ is isomorphic to a binary caterpillar. Due to the symmetry in the binary caterpillar trees, we observe the following:\\ 

\textit{The center, centroid and subtree core  coincide in binary caterpillar trees.}\\
 	
\noindent This observation shows that among all binary trees on $n$ vertices the minimum distance between any two of the above three central parts  is zero. So, it is interesting to see which trees maximize these distances among all binary trees on $n$ vertices.  There are two non-isomorphic binary trees on  $10$ vertices and both are  crg trees.  Also in any crg tree on $10$ vertices, the center, centroid and subtree core coincide. So throughout this paper, we consider binary trees on $n \geq 12$ vertices. Smith et. al. conjectured the following result in \cite{Sswy} (see Conjecture 3.10).\\

\textit{Among all binary trees on $n$ vertices, the pairwise distance between any two of center, centroid and subtree core is maximized by some trees of the family $\Omega_n$.}\\

In this paper, we  prove this conjecture.  The paper is organised in the following way: In Section \ref{Pre}, we develop some results on binary and rooted binary trees which are useful to prove our main results.  In Section \ref{Ccs}, we prove this conjecture and obtain the trees which achieve these distances.

\section{Preliminaries}\label{Pre}

Let the height of $T^n_{rg}$ be  $h\geq 1$. Then $2^h+1 \leq n \leq 2^{h+1}-1$. There exists a positive integer $\alpha$ such that $n=2^h+\alpha,$ which gives $h=\log_2(n-\alpha).$ There are two branches at the root of an rgood binary tree. The branch having maximum weight between the two, is termed as the heavier branch. If both the branches of an rgood binary tree have same weight then we say the rgood binary tree is complete. In this case any branch can be considered as heavier. The next result tells about the rooted binary trees with minimum height.
 
 \begin{lemma} \label{L-00}
Among all rooted binary trees on $n$ vertices, $ht(T^n_{rg})\leq ht(T)$ and equality holds when $T$ is a rooted binary tree in which the heights of any two pendant vertices differs by at most one.
\end{lemma}
\begin{proof}
Let $T$ be a rooted binary tree of height $h$, rooted at $r$. Suppose $T$ is not rgood. Let $v_1 \in V(T)$ such that $ht(v_1)=h=ht(T)$.

 First suppose there are no pendant vertex such that it's height differs by more than one from height of $v_1$.  Let the vertices at level $h-1$ are labled as $w_1,w_2,\ldots,w_k$ where $w_i$ is at immediate left position of $w_{i+1}$ for $i=1,2,\ldots,k-1$. Since $T$ is not rgood, the parents of the pendant vertices at level $h$ do not form a final segment at level $h-1$. Let for some $i=1,2,\ldots,k-1$, $w_i$ has no children but $w_{i+t}, t\geq 1$ has children, say $w',w''$. Delete the vertices $w',w''$ and add them as pendant vertices at $w_i$ to get a new tree $T'$. Then $ht(T')=ht(T)$. Repeat this process for $T'$ and continue till the parents of the pendant vertices at level $h$ forms a final segment at level $h-1$. Finally we get the rgood tree $T^n_{rg}$ with $ht(T^n_{rg})=h$. 
	
Now suppose there exists a pendant vertex $u_1$ such that  $ht(v_1) - ht(u_1)\geq 2$. Without loss of generality, assume that $ht(v_1)>ht(u_1)$. Let $v$ be the parent of  $v_1$. Since $T$ is a rooted binary tree, there is one more child say $v_2$ of $v$ which is a pendant vertex at level $h$. Delete the vertices $v_1,v_2$ and add them as pendant vertices at $u_1$ to get a new tree $T'$. Observe that $ht(T') \leq ht(T)$. Repeating this process we get a tree $\tilde{T}$ in which height of any two pendant vertices differs by atmost $1$ and $ht(\tilde{T})\leq ht(T')\leq ht(T)$. Now by following the argument as above, we get the rgood tree $T^n_{rg}$ with $ht(T^n_{rg})\leq ht(T)$.
\end{proof}
In the following result we determine the binary tree on $n$ vertices which has maximum diameter.
\begin{lemma} \label{L-01}
Among all binary trees on $n$ vertices the binary caterpillar  has the maximum diameter.
\end{lemma}
\begin{proof}
Let $T$ be a binary tree on $n$ vertices with diameter $k.$  Let $P:u_0u_1\ldots u_k$ be a path of maximum length in $T$. Suppose $T$ is not caterpillar. Then there exists two pendant vertices $v_1,v_2 \in V(T) - V(P)$ adjacent to $v$ such that $v$ is not on the path $P$. Delete the vertices $v_1,v_2$ and add them as a pendant vertex at $u_0$ to get a new tree $T'$. Then $diam(T')> diam(T)$. Repeat the process till a binary caterpillar  is achieved.
\end{proof}

 For  $e=\{u,v\} \in E(T)$, let $T_e(u)$ denotes the component of $T - e$ containing $u$. The following lemma is very useful .
 
 \begin{lemma}\label{Cd-03}
 Let $e=\{u,v\} \in E(T)$, then $|V(T_e(u))|>|V(T_e(v))|$ if and only if $C_d(T) \subseteq V(T_e(u)).$
 \end{lemma}
 
 \begin{proof}
Let $|V(T_e(u))|=k$ and $|V(T_e(v))|=k'$.  First suppose $|V(T_e(u))|>|V(T_e(v))|$.  Since $k>k'$, it follows that $Wt(v)=k$ and for any $w'\in V(T_e(v)), w' \neq v$, $Wt(w')>Wt(v)$. So the only possible vertex of $T_e(v)$ which may belong to $C_d(T)$ is $v$. If $Wt(u)$ is weight of the branch containing $v$, then $Wt(u)=k'<k=Wt(v)$. Otherwise, $Wt(u)$ is the weight of a branch contained in $T_e(u)$ and so  $Wt(u)\leq k-1<k=Wt(v)$. Hence $\min\{Wt(z): z\in V(T)\}\leq Wt(u)<Wt(v)$. This implies $C_d(T) \subseteq V(T_e(u))$.
 
Now suppose $C_d(T) \subseteq V(T_e(u))$. Let $w \in C_d(T)$, then $Wt(w)\geq k'$ as the branch containing $v$ at $w$ has weight at least $k'$. Since $v \not\in C_d(T)$, so  $Wt(v)>Wt(w)\geq k'$. This implies $Wt(v)$ is weight of the branch containing $u$. i.e. $Wt(v)=k$, hence $k>k'$.
\end{proof}

\begin{corollary}\label{Cor-00}
Let $v$ be the root of the rgood part of $T_{rg}^{n,l}$ and let $v'$ be the vertex in a  heavier branch of the rgood part  at $v$ such that $\{v,v'\}\in E(T_{rg}^{n,l})$. If  $l\geq \frac{n}{2}+1$  then $C_d(T_{rg}^{n,l})\subseteq \{v,v'\}$. Moreover, if the rgood part is complete then $C_d(T_{rg}^{n,l})=\{v\}.$
\end{corollary}

\begin{proof}
Let $T$ be a crg tree with $l\geq \frac{n}{2}+1$ and let $T'$ be the rgood part of $T.$ Since  $l\geq \frac{n}{2}+1$ vertices, so by Lemma \ref{Cd-03}, $C_d(T)\subseteq V(T')$.  Let $e=\{w,v\}\in E(T')$ where $w$ is not on the heavier branch of $T'$. Then a copy of $T'_e(w)$ is properly contained in $T'_e(v)$ and so $|V(T_e(v))|>|V(T_e(w))|.$ Hence by Lemma \ref{Cd-03}, $C_d(T)$ is contained in the heavier branch of $T'$. 

Suppose $e_1=\{u, v' \}\in E(T')$ where $ht(u)=2$  in  $T'$. Then a copy of $T'_{e_1}(u)$ is properly contained in $T'_{e_1}(v')$ and so $|V(T_{e_1}(v'))|>|V(T_{e_1}(u))|.$ Hence by Lemma \ref{Cd-03}, $C_d(T)\subseteq \{v, v' \}$. 

If $T'$ is complete then $T'$ has two heavier branches at $v$. Since centroid of $T$ contains either a single vertex or two adjacent vertices, so $C_d(T)=\{ v \}.$
\end{proof}

For $l\geq \frac{n}{2}+1,$ in Corollary \ref{Cor-00}, we proved that $C_d(T_{rg}^{n,l})=\{ v \}$ or $\{v' \}$ or $\{v,v'\}$. We also showed that $C_d(T_{rg}^{n,l})=\{ v \}$ if the rgood part is complete. For many values of $n$ and $l$, the other two cases will also happen. For example, it can be checked that  $C_d(T_{rg}^{12,11})=\{v'\}$ and $C_d(T_{rg}^{14,13})=\{v,v'\}$.

\begin{corollary}\label{Cor-ne}
Let $v$ be the root of the rgood part of $T_{rg}^{n,l}$ and let $v'$ be the vertex in a  heavier branch of the rgood part  at $v$ such that $\{v,v'\}\in E(T_{rg}^{n,l})$. If $n=4k$ and  $l\geq 2k+1$  then $C_d(T_{rg}^{n,l})= \{v\}$ or $\{v'\}$. 
\end{corollary}
\begin{proof}
Let $T$ be a crg tree with $n=4k$ and $l\geq 2k+1.$ By Corollary \ref{Cor-00}, $C_d(T)=\{v\}$ or $\{v'\}$ or $\{v,v'\}$. Let $e=\{v,v'\}\in E(T)$. If $C_d(T)=\{v,v'\}$ then $|V(T_e(v))|=|V(T_e(v'))|=2k$. But both $T_e(v)$ and $T_e(v')$ are binary rooted trees with roots $v$ and $v'$ , respectively and hence both must have odd number of vertices. Thus a contradiction arises, so $C_d(T)= \{v\}$ or $\{v'\}$.
\end{proof}
  
   We will now prove a result similar to Lemma \ref{Cd-03} related to  subtree core  of trees.

\begin{lemma}\label{L-06}
Let $e=\{u,v\} \in E(T)$. Then $S_c(T)\subseteq V(T_e(u))$ if and only if $f_{T_e(u)}(u)>f_{T_e(v)}(v)$.
\end{lemma}
\begin{proof} We have 
$$f_T(u)= f_{T_e(u)}(u)+f_{T_e(u)}(u)f_{T_e(v)}(v)$$ and $$f_T(v)= f_{T_e(v)}(v)+f_{T_e(u)}(u)f_{T_e(v)}(v).$$ So,  $$f_T(u)-f_T(v)=f_{T_e(u)}(u)-f_{T_e(v)}(v).$$ Now the result follows from Lemma \ref{LS-01}.
\end{proof}

\begin{corollary}\label{Cor-02}(\cite{Sswy},Proposition 1.7)
A vertex $u \in S_c(T)$ if and only if for each neighbour $v$ of $u$, $f_{T_{e_v}(u)}(u)\geq f_{T_{e_v}(v)}(v)$ where $e_v=\{v,u\}$. Furthermore if $u \in S_c(T)$ and equality holds then $v \in S_c(T)$.
 \end{corollary}
 In the following we discuss about the position of center, centroid and subtree core of $T^n_{rg}$.
 \begin{lemma} \label{L-03}
Let $v$ be the root of $T^n_{rg}$ and let $v'$ be the vertex in a  heavier branch of  $T^n_{rg}$ such that $e=\{v,v'\}\in E(T^n_{rg})$. Then center, centroid and subtree core of $T^n_{rg}$ are contained in the set $\{v,v'\}$. Moreover if $T^n_{rg}$ is complete then $C(T^n_{rg})=C_d(T^n_{rg})=S_c(T^n_{rg})=\{v\}.$
\end{lemma}
\begin{proof}
Let $P$ be a longest path of $T^n_{rg}$. Then it must go through $v$ and $C(P)=\{v\}$ or $\{v,v'\}$ depending on the length of $P$ is even or odd, respectively. So, $C(T^n_{rg})=\{v\}$ or $\{v,v'\}$.

If $n>1$, then $\deg(v)=2$. let $w'\neq v'$ and $e'=\{v,w' \}\in E(T^n_{rg})$. Since $v'$ is in the heavier branch, so $|V(T_{e'}(v))|>|V(T_{e'}(w'))|$. By Lemma \ref{Cd-03}, $C_d(T^n_{rg}) \subseteq V(T_{e'}(v))$. If $n>3$ then $\deg(v')=3.$ Let $e_1=\{v',v_1\},e_2=\{v',v_2\}\in E(T^n_{rg})$. For $i=1,2$, $T_{e_i}(v')$ contains a copy of $T_{e_i}(v_i)$. By Lemma \ref{Cd-03}, $C_d(T^n_{rg}) \subseteq V(T_{e_i}(v'))$. Hence $C_d(T^n_{rg}) \subseteq \{v, v'\}$.

Since $v'$ is in the heavier branch, so the rooted binary tree $T_{e'}(v)$ with root $v$ contains a copy of the rooted tree $T_{e'}(w')$ with root $w'$. So $f_{T_{e'}(v)}(v)>f_{T_{e'}(w')}(w')$ and hence by Lemma \ref{L-06}, $S_c(T_{rg}^n)\subseteq V(T_{e'}(v))$. Also for $i=1,2$, the rooted binary tree $T_{e_i}(v')$ with root $v'$ contains a copy of the rooted binary tree $T_{e_i}(v_i)$ with root $v_i$. So by Lemma \ref{L-06} $S_c(T_{rg}^n)\subseteq V(T_{e_i}(v'))$ for $i=1,2$. Hence $S_c(T_{rg}^n)\subseteq \{v,v'\}.$

If $T^n_{rg}$ is complete then $T^n_{rg}$ has two heavier branch at $v$ and  in this case $C(T^n_{rg})=C_d(T^n_{rg})=S_c(T^n_{rg})=\{v\}.$
\end{proof}
\begin{corollary}\label{ccsp}
Let $v$ be the root of the rgood part of $T_{rg}^{n,l}$ and let $v'$ be the vertex  on a heavier branch at $v$ with $\{v,v'\}\in E(T_{rg}^{n,l})$. Then center, centroid and subtree core of $T_{rg}^{n,l}$ lie on the path from $1$ to $v'$.
\end{corollary}

\begin{corollary}\label{Cor-03}
 Let $v$ be the root of the rgood part of $T_{rg}^{n,l}$ and let $v'$ be the vertex  on a heavier branch at $v$ with $\{v,v'\}\in E(T_{rg}^{n,l})$. Then $C(T_{rg}^{n,l})\neq \{v'\}$ . 
\end{corollary}

\begin{proof}
Let $T'$ be the rgood part of $T_{rg}^{n,l}$ and also let $ht(T')=h$.  Suppose $C({T_{rg}^{n,l}})=\{v'\}$. Then  $diam(T_{rg}^{n,l})=2(h-1)$, which  is a contradiction as  $diam(T_{rg}^{n,l})\geq 2h-1$.
\end{proof}

\subsection{Root containing subtrees}

To prove our main result, it is important to know the rooted binary trees which extremize the number of root containing subtrees. In \cite{Sswy}, the authors have obtained the rooted binary tree which maximizes the number of root containg subtrees. Here we obtain the rooted binary tree which minimizes the number of root containing subtrees.

\begin{proposition} \label{L-05}(\cite{Sswy},Corollary 3.9)
Among all rooted binary trees on $n$ vertices, $T_{rg}^n$ maximizes the number of root containing subtrees. 
\end{proposition}

For a tree $T$ with $u,v\in V(T)$, we denote the number of subtrees of $T$ containing $u$ and $v$ by $f_T(u,v).$ 
\begin{lemma}\label{L-07}
Let $T$ be a rooted binary tree with root $r$ and $x$ be a pendant vertex in $T$. Let $y$ be a vertex other than $x$ in the path joining $r$ and $x$. Then, $f_T(r,y)\geq 2f_T(r,x)$ and equality hold if and only if $y$ is adjacent to $x$.
\end{lemma}
\begin{proof}
 Suppose $x_0$ be the vertex adjacent to $x$ in $T$ and let $T_0$ be  the tree $T - x$. Then, 
\begin{equation*}
 f_T(r,x)=f_{T_0}(r,x_0)
\end{equation*} and 
\begin{equation*}
f_T(r,y)=f_{T_0}(r,y)+f_{T_0}(r,x_0)\geq 2f_{T_0}(r,x_0)=2f_T(r,x).
\end{equation*} 
The inequality holds, as any tree containing $r$ and $x_0$ must contain $r$ and $y$ and equality holds if and only if $y=x_0$.
\end{proof}
We denote the rooted binary tree on $n$ vertices with exactly two vertices at every level (except zero level) by $T_{r,2}^n.$
\begin{proposition} \label{L-08}
Among all rooted binary trees on $n$ vertices, the tree  $T_{r,2}^n$ minimizes the number of root containing subtrees.
\end{proposition}

\begin{proof}
Let $T$ be a rooted binary tree with root $r$ in which there are more than two vertices at some levels. Let $x$ be a pendant vertex of $T$ such that $ht(T)=d(r,x)$. Let $y$ be the vertex nearest to $r$ ($y$ may be same as $r$) such that every branch at $y$ contains more than two vertices. Then the path joining $r$ and $x$ must contains $y$. Let $y_0$, $y_1$ and $y_2$ be the vertices adjacent to $y.$ Let the branch at $y$ containing $y_0$ be the branch which contains $r$ (If $y=r$, then we can  take $y_1$ and $y_2$ are the only two vertices adjacent to $y$). 

Let $X$ and $Y$  be the branches  at $y$ containing $y_1$ and $y_2$, respectively and let $x$ be in the branch $Y.$  Then $X'=X-y$ is a binary rooted tree with root $y'=y_1.$  Let $T'$ be the binary rooted subtree of $T$ with root $r$, obtained by removing $X'$ from $T$ but keeping $y_1$ as a pendant vertex of it. Then $T$ can be obtained from $T'$ and $X'$ by identifying $y_1$ of $T'$ with $y'$ of $X'$. Then 
\begin{equation*}
f_T(r)=f_{T'}(r)+f_{T'}(r,y_1)(f_{X'}(y') -1).
\end{equation*}

Construct a new tree $\hat{T}$ from $T'$ and $X'$ by identifying $x$ of $T'$ with $y'$ of $X'$. Then $\hat{T}$ is a rooted binary tree with root $r$ and $|V(T)|=|V(\hat{T})|$. Then
\begin{equation*}
f_{\hat{T}}(r)=f_{T'}(r)+f_{T'}(r,x)(f_{X'}(y') -1).
\end{equation*}
So we have
\begin{equation*}
f_{T}(r)-f_{\hat{T}}(r)= (f_{X'}(y')-1)(f_{T'}(r,y_1) - f_{T'}(r,x)).
\end{equation*}
We have $f_{X'}(y')>1$ as $|V(X')|\geq 3$. Since $y_1$ is a pendant vertex so $f_{T'}(r,y_1)=f_{T'}(r,y)$. So by Lemma \ref{L-07}, $f_{T'}(r,y_1) =f_{T'}(r,y)> f_{T'}(r,x)$, hence $f_{T}(r)-f_{\hat{T}}(r)>0.$  If there are exactly two vertices at every level of $\hat{T}$ then we are done. Otherwise,  repeat the above process till we get the rooted binary tree with exactly two vertices at every level (except level zero).
\end{proof}

\begin{corollary}\label{crg}
Let $T$ be a rooted binary tree on $n$ vertices with root $r$. Then $f_T(r)\geq 3\times 2^{\frac{n-1}{2}}-2$ and equality holds if and only if $T\cong T_{r,2}^n$.
\end{corollary}
\begin{proof}
Let $r$ be the root of $T_{r,2}^n$. Suppose $u$ and $v$ are vertices adjacent to $r$ among which $u$ is pendant. Let $S_n$ be the number of subtrees of $T_{r,2}^n$ containing $r$. We have $S_1=1$ and for $n \geq 3$,
$$S_n=2S_{n-2}+2$$ where number of subtrees containing $r$ but not $v$ is $2$ and the number of subtrees containing both $r$ and $v$ is $2S_{n-2}$. We solve this recurrence relation to find the value of $S_n$. We have
\begin{align*}
S_n&=2S_{n-2}+2\\
&=2(2S_{n-4}+2)+2=2^2S_{n-4}+2+2^2\\ &\hskip .2 cm\vdots \\&= 2^{\frac{n-1}{2}}S_1+ 2+2^2+2^3+\ldots +2^{\frac{n-1}{2}}\\&=2^{\frac{n-1}{2}}+2(2^{\frac{n-1}{2}}-1)\\&= 3\times2^{\frac{n-1}{2}}-2.
\end{align*}
The result follows from Proposition \ref{L-08}.
\end{proof}
Let $r$ be the root of $T_{rg}^n$. It seems difficult to find the value of $f_{T_{rg}^n}(r)$. We will only be able to give a  bound for $f_{T_{rg}^n}(r)$ which is a solution of a nonlinear recuurence relation. Let $h$ be the height of $T_{rg}^n$ and let $m=2^{h+1}-1.$ For $n\geq 3$, $2^h-1< n \leq m$ and the rooted binary tree $T_{rg}^m$ is complete.

Let $A_h$ be the number of  subtrees  $T_{rg}^m$ containing the root $r$. We have $A_0=1$ and for $h\geq 1$, let  $u$ and $v$ be the vertices adjacent to $r$. Then $$A_h=1+A_{h-1}+A_{h-1}(1+A_{h-1})= (A_{h-1}+1)^2$$ where the first $1$ is for the subtree containing only the single vertex $r$, the second  term $A_{h-1}$ counts the number of subtrees containing $r$ and $v$ but not $u$ and the third term $A_{h-1}(1+A_{h-1})$ counts the subtrees containing $r$ and $u$. Then for $h\geq 1$, we have
$$A_{h-1}< f_{T_{rg}^n}(r)\leq A_h.$$
It will be nice to know the exact value of $f_{T_{rg}^n}(r)$.

\subsection{Solution to the recurrence $A_h=(A_{h-1}+1)^2$ for $h\geq1$, $A_0=1$}
In \cite{As}, the authors have established the solution of the following recurrence relation. \\
$x_{n+1}=x_n^2+g_n$ for $ n \geq n_0$ with boundary conditions
\begin{itemize}
\item[i)] $ x_n>0$
\item[ii)] $ |g_n|< \frac{1}{4}x_n$ \textit{and} $ x_n \geq 1$ \textit{for} $n \geq n_0$
\item[iii)] $ |\alpha_n| \geq |\alpha_{n+1}| $ \textit {for} $n\geq n_0$  
\end{itemize}
 where $\alpha_n= \ln (1+ \frac{g_n}{x_n^2}).$\\
The  authors have shown that if $x_n \in \mathbb{Z}$ and $g_n>0$, then solution of this recurrence relation is given by $x_n= \lfloor k^{2^n}\rfloor$, where $k=x_0\mbox{exp}\left({\sum_{i=o}^{\infty}(2^{-i-1}\alpha_i)} \right)$.\\
 
 We substitute $A_h+1= Y_h$ then the recurrence relation  $A_h=(A_{h-1}+1)^2$ for $h\geq1$, $A_0=1$ translates into $Y_0=2$  and $Y_{h+1}=Y_h^2+1$ for $h \geq 0$. Here $Y_h \in \mathbb{Z}$, $g_h=1>0$ and  the relation satisfies the above boundary conditions for $h_0=2$. So the solution of this is $Y_h= \lfloor k^{2^h}\rfloor$, where $k=Y_0\mbox{exp}({\sum_{i=o}^{\infty}(2^{-i-1}\alpha_i)})$ and $\alpha_i=\ln (1+ \frac{1}{Y_i^2}).$  So 
 \begin{align*}
 k&=Y_0\mbox{exp}({\sum_{i=o}^{\infty}(2^{-i-1}\alpha_i))}=2\mbox{exp}\left(\frac{1}{2}\ln(1+ \frac{1}{4})+\frac{1}{4}\ln(1+ \frac{1}{25})+\frac{1}{8}\ln(1+ \frac{1}{676})+\cdots \right)\\&=2\mbox{exp}\left(\frac{1}{2}\ln( \frac{5}{4})+\frac{1}{4}\ln( \frac{26}{25})+\frac{1}{8}\ln( \frac{677}{676})+\cdots\right)=2.25851845\cdots .
\end{align*}

 Hence $A_h=\lfloor k^{2^h}\rfloor -1$ where $k=2.25851845\cdots $.
 
\section{Center, Centroid and Subtree core}\label{Ccs}

In this section we obtain the binary trees which maximize the pairwise distances between the central parts center, centroid and subtree core over all binary trees on $n$ vertices. We first consider the pair center and centroid.

\subsection{Center and centroid}

\begin{theorem}\label{th-00}
Among all binary trees on $n$ vertices, the distance between  center and  centroid is maximized by a crg tree.
\end{theorem}
\begin{proof}
Let $T$ be a binary tree on $n$ vertices with $d_T(C,C_d)\geq 1.$ Let  $u\in C(T)$ and $v\in C_d(T)$ such that $d_T(C,C_d)=d(u,v)$. Let $e=\{v,w\} \in E(T)$ such that $w$ lies on the path joining $u$ and $v$. Let $|V(T_e(v))|=k$. The component $T_e(v)$ is a rooted tree with root $v$. Since $C_d(T)\subseteq T_e(v)$ so by Lemma \ref{Cd-03}, $|V(T_e(v))|>|V(T_e(w))|$.

If $T_e(v)$ is a rgood binary tree then rename the tree $T$ as $T'$. Otherwise, form a new tree $T'$ from $T$ by replacing the component $T_{e}(v)$ with $T^k_{rg}$ rooted at $v$. Since $|V(T'_e(w))|<k=|V(T'_e(v))|=|V(T^k_{rg})|$, so  by Lemma \ref{Cd-03}, $C_d(T')\subseteq V(T'_e(v))$. Also by Lemma \ref{L-00}, $ht(T_e(v))\geq ht(T^k_{rg})$. So  $C(T')$ either  same as $C(T)$ or moves away from the vertex $v$. Hence, $d_{T'}(C,C_d) \geq d_T(C,C_d)$.
	
If $T' \in \Omega_n$ then the result follows. Otherwise let $|V(T'_{e}(w))|=l$.  Construct a new tree $T''$ from $T'$ by replacing $T'_{e}(w)$  with $T_{r,2}^l$ rooted at  $w$. Observe that $T'' \in \Omega_n$. In $T''$ the length of the longest path is more than  the length of the longest path of $T'$ and the increment occurs in a branch at $w$ containing the center. So, $C(T'')$ moves away from $v$. Also, $|V(T''_{e}(v))|> |V(T''_{e}(w))|$ and $T''_{e}(v)$ is same as $T'_{e}(v)$. So $C_d(T'')=C_d(T')$. Hence $d_{T''}(C,C_d) \geq d_{T'}(C,C_d) \geq d_T(C,C_d)$. This proves the result.
\end{proof}
\begin{theorem}
Among all crg trees on $n\geq 12$ vertices, the distance between center and centroid is maximized by the  tree $T_{rg}^{n,l},$ where $l=2\lceil\frac{n}{4}\rceil+1$.
\end{theorem}
\begin{proof}
Let $v$ be the root of the rgood part of $T_{rg}^{n,l}$ and let $d_{T_{rg}^{n,l}}(C,C_d)=\alpha$. We consider two cases depending on whether $n$ is of the form $4k$ or $4k+2$.

\textbf{Case I:} $n=4k$ for some $k\geq 3.$\\
In this case $l=2k+1.$ In $T_{rg}^{n,l}$, the rgood part  has $2k+1$ vertices and the caterpilar part has $2k$ vertices. The weight $Wt(v)=2k-1$ and the weight of any other vertex of $T_{rg}^{n,l}$ is greater than $2k-1$. Following the numbering of vertices mentioned in Section \ref{num}, we have $C_d(T_{rg}^{n,l})=\{v\}=\{k+1\}$. The diameter of the caterpilar part is $k$ and the height of the the rgood part is less than $k$. So, $C(T_{rg}^{n,l})$ lies in the path from $1$ to $k+1$.

First consider the  trees $T_{rg}^{n,l},T_{rg}^{n,l-2},\ldots,T_{rg}^{n,5}$. Note that  $T_{rg}^{n,5}$ is a binary caterpilar. Then $C_d(T_{rg}^{n,l-i})=\{k+1\}$ for $0\leq i\leq l-5.$ In the above sequence of  trees,  the center lies in the path from $1$ to $k+1$.  In $T_{rg}^{n,l-i}$, if the vertex numbered $u$ is the central vertex nearest to $k+1$, then the central vertex in $T_{rg}^{n,l-i-2}$ nearest to $k+1$ is either $u$ or $u+1$. So, $d_{T_{rg}^{n,l-i}}(C,C_d)\geq d_{T_{rg}^{n,l-i-2}}(C,C_d)$ for $0\leq i\leq l-7.$  

Now consider the sequence of trees $T_{rg}^{n,l},T_{rg}^{n,l+2},\ldots,T_{rg}^{n,n-1}.$  For $0\leq j \leq n-l-1,$ let $v_j$ be the root of the rgood part of  $T_{rg}^{n,l+j}$ and $v_j'$ be the vertex in a heavier branch of the rgood part of $T_{rg}^{n,l+j}$ adjacent to $v_j.$  If $d_{T_{rg}^{n,l}}(C,C_d)=\alpha=0$ then $C(T_{rg}^{n,l})=\{k+1\}$ or $\{k,k+1\}$. So $\alpha=0$ implies $n\leq 16$  and in these cases it can be checked that $d_{T_{rg}^{n,l+j}}(C,C_d)=0$ for $0\leq j \leq n-l-1.$
 If  $\alpha\geq 1$ then $n\geq 20$. If $j'$ is the smallest positive integer such that $C_d(T_{rg}^{n,l+j'})=\{v_{j'}\}$ then $C_d(T_{rg}^{n,l+j'+p})=\{v_{j'}\}$ for $0\leq p \leq n-l-j'-1$. Since $n\geq 20$, so $C_d(T_{rg}^{n,l+2})=\{v_2\}$. Then $d_{T_{rg}^{n,l+2}}(C,C_d)=\alpha$ or $\alpha-1$. If $d_{T_{rg}^{n,l+2}}(C,C_d)=\alpha$ then $d_{T_{rg}^{n,l+4}}(C,C_d)=\alpha$ or $\alpha-1$. If $d_{T_{rg}^{n,l+4}}(C,C_d)=\alpha$ then $C_d(T_{rg}^{n,l+4})=\{v_4'\}$ and $C_d(T_{rg}^{n,l+j})=v_j'$ for  $4\leq j \leq n-l-1$. Hence the distance between center and centroid  of all the trees in the above sequnce is at most $\alpha.$ This prove the result for the case $n=4k.$

\textbf{Case II:} $n=4k+2$ for some $k\geq 3.$\\
A similar argument can be given to prove this case. This completes the proof.
\end{proof}

We will now find the distance between center and centroid of $T_{rg}^{n,l}$, for $l=2\lceil\frac{n}{4}\rceil+1$. Let $h$ be the height of the rgood part of $T_{rg}^{n,l}$. Then $h$ is the smallest positive integer such that $l\leq 2^{h+1}-1.$ This implies $\lceil\frac{n}{4}\rceil\leq  2^h -1.$

The caterplilar part of  $T_{rg}^{n,l}$ contains $2\lfloor \frac{n}{4}\rfloor$ vertices. So, the diameter of the caterpilar part is $\lfloor \frac{n}{4}\rfloor$ and hence the root of the rgood part is numbered by  $\lfloor \frac{n}{4}\rfloor+1$. The central vertex which is nearest to the root of the  rgood part is numbered by  $\left\lfloor\frac{\lfloor \frac{n}{4}\rfloor+1+h}{2}\right\rfloor+1$. Thus the distance between center and centroid of $T_{rg}^{n,l}$ is $\lfloor \frac{n}{4}\rfloor - \left\lfloor\frac{\lfloor \frac{n}{4}\rfloor+1+h}{2}\right\rfloor$, where $h$ is the smallest positive integer such that $\lceil\frac{n}{4}\rceil\leq  2^h -1.$ This leads to the following corollary.
\begin{corollary}
Let $T$ be a binary tree on $n$ vertices and let $h$ be the smallest positive integer such that $\lceil\frac{n}{4}\rceil\leq  2^h -1$. Then 
$$d_T(C,C_d)\leq \left\lfloor \frac{n}{4}\right\rfloor - \left\lfloor\frac{\lfloor \frac{n}{4}\rfloor+1+h}{2}\right\rfloor$$ 
and equality happens if $T \cong T_{rg}^{n,l}$ where $l=2\lceil\frac{n}{4}\rceil+1$.
\end{corollary}

\subsection{Center and Subtree core}

\begin{theorem}\label{th-01}
Among all binary trees on $n$ vertices, the distance between  center and  subtree core is maximized by a crg tree.
\end{theorem}

\begin{proof}
Let $T$ be a binary tree on $n$ vertices with $d_T(C,S_c)\geq 1.$ Let  $u\in C(T)$ and $v\in S_c(T)$ such that $d_T(C,S_c)=d(u,v)$. Let $e=\{v,w\} \in E(T)$ such that $w$ lies on the path joining $u$ and $v$. Let $|V(T_e(v))|=k$. The component $T_e(v)$ is a rooted binary tree with root $v$. Since $S_c(T)\subseteq V(T_e(v))$ so by Lemma \ref{L-06}, $f_{T_e(v)}(v)>f_{T_e(w)}(w)$.

If $T_e(v)$ is a rgood binary tree then rename the tree $T$ by $T'$. Otherwise, form a new tree $T'$ from $T$ by replacing the component $T_{e}(v)$ with $T^k_{rg}$ rooted at $v$. By Propositon \ref{L-05}, $f_{T^k_{rg}}(v)\geq f_{T_e(v)}(v)>f_{T_e(w)}(w)$ and hence by Lemma \ref{L-06}, $S_c(T')\subseteq V(T'_e(v)).$ We have by Lemma \ref{L-00}, $ht(T_e(v))\geq ht(T^k_{rg})$. So  $C(T')$ either  same as $C(T)$ or moves away from the vertex $v$. Hence, $d_{T'}(C,S_c) \geq d_T(C,S_c)$.
	
If $T' \in \Omega_n$ then the result follows. Otherwise let $|V(T'_{e}(w))|=l$.  Construct a new tree $T''$ from $T'$ by replacing $T'_{e}(w)$  with $T_{r,2}^l$ rooted at  $w$. Observe that $T'' \in \Omega_n$. In $T''$ the length of the longest path is more than  the length of the longest path of $T'$ and the increment occurs in a branch at $w$ containing the center. So, $C(T'')$ moves away from $v$.  By Proposition \ref{L-08}, $f_{T_{r,2}^l}(w)< f_{T'_{e}}(w)$.
So $f_{T''_{e}(v)}(v)>f_{T''_{e}(w)}(w)$ and hence by Lemma \ref{L-06}, $S_c(T'')\subseteq V(T''_{e}(v))$.  Thus $d_{T''}(C,S_c) \geq d_T(C,S_c)$. This proves the result.
\end{proof}
 
\begin{theorem}\label{th-m}
In  any crg tree $T_{rg}^{n,l}$, the centroid lies in the path connecting the center and the subtree core.
\end{theorem}
 \begin{proof}
 In a binary caterpilar tree on $n$ vertices the center, centroid and subtree core are same. So we can consider crg trees which are not caterpilar. Let $T$ be a crg non-caterpilar tree, and let $T'$ be the rgood part of $T$. Let $v$ be the root of $T'$ and let $v'$ be the vertex in a heavier branch of $T'$ such that $\{v,v'\}\in E(T')$. By Corollary \ref{ccsp},  the center, centroid and subtree core of $T$ lie in the path from $1$ to $v'$.
 
 Let $w$ be the centroid vertex of $T$ nearest to $v'$ ($w$ may be same as $v'$). Let $w'$ be the vertex adjacent to $w$ and lies in the path from $1$ to $w$. Let $e=\{w', w\}\in E(T)$. Then by Lemma \ref{Cd-03}, $|V(T_{e}(w)|\geq |V(T_{e}(w')|$. If $w=v$ or $v'$ then $T_{e}(w)$ is a rgood binary tree with root $w$. By Propositon \ref{L-05}, $f_{T_{e}(w)}(w) >f_{T_{e}(w')}(w')$ and hence by Lemma \ref{L-06}, $S_c(T)\subseteq V(T_e(w))$. If $w$ is neither $v$ nor $v'$ then $T_{e}(w')$ is a rooted binary tree in which every level has exactly two vertices( except the zero level). By Propositon \ref{L-08}, $f_{T_{e}(w)}(w) >f_{T_{e}(w')}(w')$ and hence by Lemma \ref{L-06}, $S_c(T)\subseteq V(T_e(w))$. Hence $S_c(T)$ lies in the path between $C_d(T)$ and $v'$.
 
Let $u$ be the central vertex of $T$ nearest to the vertex $1$. Let $u'$ be the vertex  lies in the path from $1$ to $u$ ($u'$ may be same as $1$) such that $e_1=\{u',u\}\in E(T)$. By Corollary \ref{Cor-03}, $u$ lies in the path from $1$ to $v$ and the component $T_{e_1}(u')$ is a rooted binary tree with root $u'$  in which every level has exactly two vertices (except the zero level). Hence $|V(T_{e_1}(u')|< |V(T_{e_1}(u)|$ and by Lemma \ref{Cd-03}, $C_d(T)\subseteq V(T_{e_1}(u).$ Thus $C_d(T)$ lies  in the path from $C(T)$ to $v'$. This completes the proof.
 \end{proof}
 
 \begin{corollary}\label{close}
 In a crg tree $T_{rg}^{n,l}$, among center, centroid and subtree core, the  center is nearest to the vertex $1$. 
 \end{corollary}
 \begin{proof}
 We rename the crg tree $T_{rg}^{n,l}$ as $T$. Let $u\in C(T)$ and $v \in S_c(T)$ such that $d(1,C(T))=d(1,u)$ and $d(1,S_c)=d(1,v)$. We show that $d(1,v)\geq d(1,u)$. Suppose $d(1,v)<d(1,u)$. Let $w$ be the vertex adjacent to $u$ in the path joining $1$ and $u$. Consider the edge $e=\{w,u\}$. Then $v \in V(T_e(w))$. Also, as $u \in C(T)$ is the central vertex nearest to $1$, $k=|V(T_e(w))|<|V(T_e(u))|$. Note that the tree $T_e(w)$ is $T_{r,2}^k$. So by Lemma \ref{L-08}, $f_{T_e(w)}(w)<f_{T_e(u)}(u).$ Hence by Lemma \ref{L-06}, $v \in V(T_e(u))$, which is a contradiction.
 \end{proof}

 Let $l$ be an odd integer and let $r$ be the root of $T_{rg}^{l}$. We denote the number $f_{T_{rg}^l}(r)$ by $R_l$.
 \begin{lemma}\label{msc}
 Let $n\geq 12$ be even and let $l$ be the smallest positive odd number such that $R_l>3\times2^{\frac{n-l-1}{2}} -2.$ For an even integer $i$ with $2\leq i \leq l-3$, let $e=\{\frac{n-l+3}{2},\frac{n-l+5}{2}\} \in E(T_i)$ where $T_i=T_{rg}^{n,l-i}.$ Then $S_c(T_i)\subseteq V(T_{i_e}(\frac{n-l+3}{2}))$ or $S_c(T_i)=\{\frac{n-l+3}{2},\frac{n-l+5}{2}\}.$
 \end{lemma}
 \begin{proof}
 Since $n\geq 12$ so $l\geq 5$ and $T_i$ is defined for every even integer $i$ with $2\leq i \leq l-3.$ Also since $l$ is the smallest positive odd integer such that $R_l>3.2^{\frac{n-l-1}{2}} -2$ so $R_{l-i}\leq R_{l-2}\leq 3\times 2^{\frac{n-l+1}{2}} -2.$
 
The component $T_{i_e}(\frac{n-l+5}{2})$ of $T_i-e$ is a binary rooted tree on $l-2$ vertices with root $\frac{n-l+5}{2}$ ($T_{i_e}(\frac{n-l+5}{2})$ is a rgood tree if $i=2$). By Proposition \ref {L-05}, $ f_{T_{i_e}}(\frac{n-l+5}{2}) \leq R_{l-2}$ for $2 \leq i \leq l-5$. The component $T_{i_e}(\frac{n-l+3}{2})$ of $T_i-e$ is a binary rooted tree on $n-l+2$ vertices with root $\frac{n-l+3}{2}.$ By Corollary \ref{crg}, $f_{T_{i_e}}(\frac{n-l+3}{2})=3\times 2^\frac{n-l+1}{2}-2.$

Thus we have  $$f_{T_{i_e}}(\frac{n-l+3}{2})=3.2^\frac{n-l+1}{2}-2\geq R_{l-2}\geq f_{T_{i_e}}(\frac{n-l+5}{2})$$ for $2 \leq i \leq l-5$.  If $S_c(T_i)\neq \{\frac{n-l+3}{2},\frac{n-l+5}{2}\}$ then $f_{T_{i_e}}(\frac{n-l+3}{2})>f_{T_{i_e}}(\frac{n-l+5}{2})$. Then by Lemma \ref{L-06}, $S_c(T_i)\subseteq V(T_{i_e}(\frac{n-l+3}{2}))$. This comletes the proof.
 \end{proof}
\begin{lemma}\label{msc1}
Let $n\geq 12$ be even and let $l$ be the smallest positive odd number such that $R_l>3\times2^{\frac{n-l-1}{2}} -2.$
Then  $S_c(T_{rg}^{n,l+j})=\{v_j\}$ where $v_j$ is the root of the rgood part of $T_{rg}^{n,l+j}$ for $0\leq j\leq 2$.
\end{lemma} 
\begin{proof}
Since $T_{rg}^{n,5}$ is a binary caterpillar and $n\geq 12$ so $l\geq 7.$ For $0\leq j\leq 2$, let $v_j'$ be the vertex in a heavier branch of of the rgood part of $T_{rg}^{n,l+j}$ with $e_j=\{v_j,v_j' \}\in E(T_{rg}^{n,l+j}).$ We have $S_c(T_{rg}^{n,l+j})\subseteq \{v_j,v_j' \}$ as $R_l>3\times2^{\frac{n-l-1}{2}} -2.$ The component $T_{e_j}(v_j')$ of $T_{rg}^{n,l+j} - e_j$ is a rgood tree with root $v_j'$ and has at most $l-2$ vertices. The component $T_{e_j}(v_j)$ of $T_{rg}^{n,l+j} - e_j$ is a rooted binary tree with root $v_j$ and has at least $n-l+2$ vertices. Also at level one of $T_{e_j}(v_j)$ more than two vertices. Thus we have,
$$f_{T_{e_j}(v_j')}(v_j')\leq R_{l-2}\leq 3\times2^{\frac{n-l+1}{2}} -2 < f_{T_{e_j}(v_j)}(v_j)$$ for $0\leq j\leq 2$. Hence by Lemma \ref{L-06},  $S_c(T_{rg}^{n,l+j})=\{v_j\}$ for $0\leq j\leq 2$.
\end{proof}
 \begin{theorem}
 Let $n\geq 12$ be even and let $l$ be the smallest positive odd number such that $R_l>3\times 2^{\frac{n-l-1}{2}} -2.$ Among all crg trees on $n$ vertices, the distance between center and subtree core is maximized by the  tree $T_{rg}^{n,l}.$ 
 \end{theorem}
\begin{proof}
Let $v$ be the root of the rgood part of  $T_{rg}^{n,l}.$ Then by Lemma \ref{msc1}, $S_c(T_{rg}^{n,l})=\{v\}.$  Following the numbering of vertices mention in Section \ref{num},  $v$ is numbered as $\frac{n-l+3}{2}$ in $T_{rg}^{n,l}$. Let the vertex numbered  as $u$ be the central vertex of $T_{rg}^{n,l}$ nearest to the vertex $\frac{n-l+3}{2}$. Then by Corollary \ref{close}, $u$ lies on the path joining $1$ and $\frac{n-l+3}{2}$ and $$d_{T_{rg}^{n,l}}\left(u,\frac{n-l+3}{2}\right)= d_{T_{rg}^{n,l}}(C,S_c).$$

Let $i$ be an even integer with $2\leq i \leq l-3.$ Consider the crg tree $T_{rg}^{n,l-i}$. Since the center of a tree is same as the center of every longest path in it, so  the central vertex of $T_{rg}^{n,l-i}$ nearest to $\frac{n-l+3}{2}$ is $u+k$ for some $k\geq 0.$  Also by Lemma \ref{msc}, $S_c(T_{rg}^{n,l-i})=\{\frac{n-l+3}{2},\frac{n-l+5}{2}\}$ or lies on the path from $1$ to $\frac{n-l+3}{2}$. Hence, we have

$$d_{T_{rg}^{n,l-i}}(C,S_c)\leq d_{T_{rg}^{n,l-i}}(u,S_c)\leq d_{T_{rg}^{n,l-i}}\left(u,\frac{n-l+3}{2}\right)= d_{T_{rg}^{n,l}}(C,S_c)$$ for $2\leq i \leq l-3.$

Consider the sequence of trees $T_{rg}^{n,l+j}$ for $0\leq j\leq n-l-1$ with $j$ even. Then by Lemma \ref{L-06} and Corollary \ref{ccsp}, $S_c(T_{rg}^{n,l+j})\subseteq \{w,w'\}$ for $0\leq j\leq n-l-1$ where $w$ is the root of the rgood part of  $T_{rg}^{n,l+j}$ and  $w'$ is the vertex in a heavier branch with $e=\{w,w'\}\in E(T_{rg}^{n,l+j}).$ In $T_{rg}^{n,l+j}$, $w$ is numbered as $\frac{n-l-j+3}{2}.$ Let $d_{T_{rg}^{n,l}}(C,S_c)=\alpha.$ We have two cases :

\textbf{Case I:} $\alpha\geq 1$\\
By Lemma \ref{msc1}, $d_{T_{rg}^{n,l+2}}(C,S_c)=\alpha$ or $\alpha-1$. Let $j'$ be the smallest positive even integer such that $d_{T_{rg}^{n,l+j'}}(C,S_c)=0.$  Then $d_{T_{rg}^{n,l+k}}(C,S_c)=0 \; \mbox{or} \; 1$, for $j'\leq k \leq n-l-1$  and $d_{T_{rg}^{n,l+k}}(C_d,S_c)\leq \alpha$ for $0\leq k \leq j'-2.$ Hence 
$$d_{T_{rg}^{n,l+j}}(C,S_c)\leq  d_{T_{rg}^{n,l}}(C,S_c)$$ for $2\leq j \leq n-l-1.$\\

\textbf{Case II:} $\alpha=0$\\
Since $n$ is even, so $n=4k $ or $4k+2$ for some $k.$ So $l\leq 2k+1$ as $l$ is the smallest positive odd number such that $R_l>3\times 2^{\frac{n-l-1}{2}} -2.$ It can be checked that $C(T_{rg}^{14,7})=\{4\}$ and $S_c(T_{rg}^{14,7})=\{5\}$. So $d_{T_{rg}^{14,7}}(C,S_c)=1$ and hence for $n=4k+2$ for some $k\geq 3,$ $d_{T_{rg}^{n,l}}(C,S_c)\geq 1.$
It can also be checked that $C(T_{rg}^{20,11})=\{5\}$ and $S_c(T_{rg}^{20,11})=\{6\}$. So $d_{T_{rg}^{20,11}}(C,S_c)=1$ and hence for $n=4k$ for some $k\geq 5,$ $d_{T_{rg}^{n,l}}(C,S_c)\geq 1.$

If $n=12$ then $l=7$ and it can be easily checked that $d_{T_{rg}^{12,7}}(C,S_c)=d_{T_{rg}^{12,9}}(C,S_c)=d_{T_{rg}^{12,11}}(C,S_c)=0.$ If $n=16$ then $l=9$ and it also can be cehecked that $d_{T_{rg}^{16,9}}(C,S_c)=d_{T_{rg}^{16,11}}(C,S_c)=d_{T_{rg}^{16,13}}(C,S_c)= d_{T_{rg}^{16,15}}(C,S_c)=0.$ Hence if $d_{T_{rg}^{n,l}}(C,S_c)=0$ then $d_{T_{rg}^{n,l+j}}(C,S_c)=0$ for $2\leq j \leq n-l-1.$ This completes the proof.
\end{proof} 
\begin{corollary}
Let $T$ be a binay tree on $n\geq 12$ vertices. Let $r$ be the root of the rooted binary tree $T_{rg}^{l}$ and let $l$ be the smallest postive integer such that  $f_{T_{rg}^{l}}(r)>3\times 2^{\frac{n-l-1}{2}} -2.$ Then $$d_T(C,S_c)\leq d_{T_{rg}^{n,l}}(C,S_c).$$ 
\end{corollary} 

 \subsection{Centroid and Subtree core}
 
\begin{theorem}\label{lt1}
Among all binary trees on $n$ vertices, the distance between centroid and subtree core is maximized by a crg tree.
\end{theorem}
 \begin{proof}
 Consider a binary tree $T$ on $n$ vertices with $d_T(C_d,S_c)\geq 1$. Our aim is to construct a crg tree $\tilde{T} \in \Omega_n$ such that $d_{\tilde{T}}(C_d,S_c)\geq d_T(C_d,S_c)$. Let $u\in C_d(T)$ and $v \in S_c(T)$ such that $d_T(C_d,S_c)=d(u,v)$.  Let $u'$ and $v'$ be the vertices adjacent to $u$ and $v$ respectively, and lie on the path joining $u$ and $v$. Let $e_1=\{u,u'\}, e_2=\{v',v\} \in E(T)$.
 
 Let $|V(T_{e_2}(v))|=k$. The component $T_{e_2}(v)$ is a rooted binary tree with root $v$. Since $S_c(T)\subseteq V(T_{e_2}(v))$ so by Lemma \ref{L-06}, $f_{T_{e_2}(v)}(v)>f_{T_{e_2}(v')}(v')$. If $T_{e_2}(v)$ is an rgood binary tree then rename the tree $T$ by $T'$. Otherwise, form a new tree $T'$ from $T$ by replacing the component $T_{e_2}(v)$ with $T^k_{rg}$ rooted at $v$. By Propositon \ref{L-05}, $f_{T^k_{rg}}(v)\geq f_{T_{e_2}(v)}(v)>f_{T_{e_2}(v')}(v')$ and hence by Lemma \ref{L-06}, $S_c(T')\subseteq V(T'_{e_2}(v)).$ Since $C_d(T)\subseteq V(T_{e_1}(u))$ so by Lemma \ref{Cd-03}, $|V(T_{e_1}(u))| >|V(T_{e_1}(u'))|$. As $d_T(C_d,S_c)\geq 1$ and $V(T)=V(T')$ so $C_d(T)=C_d(T')$. Hence, $d_{T'}(C_d,S_c) \geq d_T(C_d,S_c)$.
 
If $T' \in \Omega_n$ then the result follows. Otherwise let $|V(T'_{e_2}(v'))|=l$.  Construct a new tree $T''$ from $T'$ by replacing $T'_{e_2}(v')$  with $T_{r,2}^l$ rooted at  $v'$. Observe that $T'' \in \Omega_n$.  By Proposition \ref{L-08}, $f_{T_{r,2}^l}(v')< f_{T'_{e_2}}(v')$. So $f_{T''_{e_2}(v)}(v)>f_{T''_{e_2}(v')}(v')$ and hence by Lemma \ref{L-06}, $S_c(T'')\subseteq V(T''_{e_2}(v))$. Also we can construct $T''$ from $T'$ step wise such that in each step the centroid is same as the centroid of $T'$ or moves away from $v.$ For that choose a longest path $P$ starting from $v'$ containg the centroid of $T$ in the binary rooted tree $T'_{e_2}(v')$ with root $v'$. Let $x$ be the end poind of the path $P$. Delete two pendant vertices from same parents, where the parent is not on the path $P$ and add them as a pendant vertices at $x.$ Continue this process till $T'_{e_2}(v')$ becomes the tree $T_{r,2}^l$ and we reach the tree $T''$. In each step of the process, the centroid is either same as the centroid of the tree in the previous step or moves away from $v.$ Hence, $d_{T''}(C_d,S_c)\geq d_T(C_d,S_c)$. This completes the proof.
 \end{proof}
 
 \begin{theorem}
 Let $T$ be a binay tree on $n\geq 12$ vertices. Let $r$ be the root of the rooted binary tree $T_{rg}^{l}$ and let $l$ be the smallest postive integer such that  $f_{T_{rg}^{l}}(r)>3\times 2^{\frac{n-l-1}{2}} -2.$ Then 
 \begin{equation*}
d_T(C_d,S_c)\leq\begin{cases}
 d_{T_{rg}^{n,l}}(C_d,S_c)  & \; \mbox{if}\;\;  d_{T_{rg}^{n,l}}(C_d,S_c)\geq 1\\
 1   & \textit{otherwise}. 
\end{cases}
\end{equation*}
 \end{theorem}
\begin{proof}
Since $n$ is even, so $n=4k $ or $4k+2$ for some $k.$ So $l\leq 2k+1$ as $l$ is the smallest positive odd number such that $R_l>3\times 2^{\frac{n-l-1}{2}} -2.$ Then $C_d(T_{rg}^{n,l})$ lies in the path from $1$ to $\frac{n-l+3}{2}$. Let the vertex numbered $u$ be the centroid vertex of $T_{rg}^{n,l}$ nearest to the vertex $\frac{n-l+3}{2}$. Then by  Lemma \ref{msc1}, $$d_{T_{rg}^{n,l}}\left(u,\frac{n-l+3}{2}\right)= d_{T_{rg}^{n,l}}(C_d,S_c).$$

Let $i$ be an even integer with $2\leq i \leq l-3.$ Consider the crg tree $T_{rg}^{n,l-i}$.   Then the centroid vertex of $T_{rg}^{n,l-i}$ nearest to $\frac{n-l+3}{2}$ is $u$ for $2\leq i \leq l-3.$ Also by Lemma \ref{msc}, $S_c(T_{rg}^{n,l-i})=\{\frac{n-l+3}{2},\frac{n-l+5}{2}\}$ or lies on the path from $1$ to $\frac{n-l+3}{2}$. So we have

$$d_{T_{rg}^{n,l-i}}(C_d,S_c)= d_{T_{rg}^{n,l-i}}(u,S_c)\leq d_{T_{rg}^{n,l-i}}\left(u,\frac{n-l+3}{2}\right)= d_{T_{rg}^{n,l}}(C_d,S_c)$$ for $2\leq i \leq l-3.$

Consider the sequence of trees $T_{rg}^{n,l+j}$ for $0\leq j\leq n-l-1$ with $j$ even. Then by Lemma \ref{L-06} and Corollary \ref{ccsp}, $S_c(T_{rg}^{n,l+j})\subseteq \{w,w'\}$ for $0\leq j\leq n-l-1$ where $w$ is the root of the rgood part of  $T_{rg}^{n,l+j}$ and  $w'$ is the vertex in a heavier branch with $e=\{w,w'\}\in E(T_{rg}^{n,l+j}).$ In $T_{rg}^{n,l+j}$, $w$ is numbered as $\frac{n-l-j+3}{2}.$ Let $d_{T_{rg}^{n,l}}(C_d,S_c)=\alpha.$ We have two cases :

\textbf{Case I:} $\alpha\geq 1$\\
 Let $j'$ be the smallest positive even integer such that $d_{T_{rg}^{n,l+j'}}(C_d,S_c)=0.$ Then $d_{T_{rg}^{n,l+k}}(C_d,S_c)=0 \; \mbox{or} \; 1$, for $j'\leq k \leq n-l-1$  and $d_{T_{rg}^{n,l+k}}(C_d,S_c)\leq \alpha$ for $0\leq k \leq j'-2.$ Hence 
$$d_{T_{rg}^{n,l+j}}(C_d,S_c)\leq  d_{T_{rg}^{n,l}}(C_d,S_c)$$ for $2\leq j \leq n-l-1.$\\

\textbf{Case II:} $\alpha=0$\\
In this case $d_{T_{rg}^{n,l+j}}(C_d,S_c)=0 \; \mbox{or} \; 1$,  for $2\leq j \leq n-l-1.$\\

Hence the result follows from Theorem \ref{lt1}.
\end{proof}

\noindent{\bf Addresses}:\\

\noindent 1) School of Mathematical Sciences,\\
National Institute of Science Education and Research (NISER), Bhubaneswar,\\
P.O.- Jatni, District- Khurda, Odisha - 752050, India\medskip

\noindent 2) Homi Bhabha National Institute (HBNI),\\
Training School Complex, Anushakti Nagar,\\
Mumbai - 400094, India\medskip

\noindent E-mails: dinesh.pandey@niser.ac.in, klpatra@niser.ac.in

\end{document}